\documentclass[krantz1,ChapterTOCs]{krantz} 
\usepackage{fixltx2e,fix-cm}
\usepackage{amssymb}
\usepackage{amsmath}
\usepackage{graphicx}
\usepackage{makeidx}
\usepackage{multicol}
\usepackage{lipsum}
\usepackage{pinlabel}

\newcommand{\showhide}[1]{#1} 
\showhide{\usepackage{tikz}}

\AtBeginDocument{\renewcommand{\bibname}{Resource List}}
\HeadingsChapterSection

\makeindex  

\newtheorem{thm}{Theorem}[section] 
\newtheorem{lem}[thm]{Lemma}
\newtheorem{cor}[thm]{Corollary}

\newtheorem{prop}[thm]{Proposition} 

\newtheorem{defn}[thm]{Definition}

\newtheorem{ex}[thm]{Example}

\numberwithin{equation}{section}

\numberwithin{figure}{section}

\newcommand{\Ext}{\ensuremath{\mathsf \Lambda}}
\newcommand{\nbc}{\ensuremath{\text{nbc}}} 

\newcommand{\OS}{\ensuremath{\mathsf A}}

\newcommand{\R}{\ensuremath{\mathsf R}}

\renewcommand{\theenumi}{\roman{enumi}}

\renewcommand{\k}{\ensuremath{\Bbbk}}
\renewcommand{\bar}{\overline}

\newcommand{\cl}{\operatorname{cl}}

\begin{document}
\frontmatter

\chapter*{Symbols} 
\begin{symbollist}{0000000}
\symbolentry{$T(G;x,y)$}{The Tutte polynomial of  $G$}


\symbolentry{$\mathcal{G}(M)$}{The $\mathcal{G}$-invariant of  $M$}

\symbolentry{$\mathsf{A}(M)$}{The Orlik-Solomon algebra of $M$}

\end{symbollist}

\mainmatter

\chapterauthor{Michael J. Falk}{Department of Mathematics and Statistics, \\ 
Northern Arizona University, \\ 
Flagstaff, Arizona 86011-5717, USA  \\
\texttt{michael.falk@nau.edu}\\}

\chapterauthor{Joseph P.S. Kung}{Department of Mathematics \\  University of North Texas\\
Denton, Texas 76203, USA, \\
\texttt{kung@unt.edu}\\}

\chapter{Algebras and valuations related to the Tutte polynomial}

\section{Synopsis}

We give brief introduction to three topics in which Tutte or characteristic polynomials play a role:   
\begin{enumerate} 
\item[$\bullet$]  The Orlik-Solomon algebra of a simple matroid 
\item[$\bullet$]  The  $\mathcal{G}$-invariant and other valuative invariants on matroid polytopes  
\item[$\bullet$]  Coalgebras constructed from graphs and matroids.  
\end{enumerate} 
The three topics are essentially independent and can be read separately.

\section{Introduction}  

In addition to playing a central role in graph and matroid theory, Tutte polynomials or specializations have also appeared in other areas.   We shall discuss briefly three of these areas.  We begin with a description of the Orlik-Solomon algebra of a simple matroid $M.$   
These algebras are quotients of exterior algebras by ideals generated by relations defined using the circuits of $M$ and the dimension of the subspace of elements of a given grade is a coefficient of the characteristic polynomial.  We continue with an exposition of valuative invariants on polytopes defined by bases of a matroid.  The Tutte polynomial is  
a valuative invariant and in particular, it is a specialization of a universal valuative invariant called the 
$\mathcal{G}$-invariant.  This is a new but important area and research is ongoing; thus our exposition can only be tentative.  We end with a short account of coalgebras associated with graphs and matroids. 

Throughout this chapter, we will usually abbreviate the set $\{a,b,\ldots,d\}$ by $ab\ldots d;$ for example, $123$ is the set $\{1,2,3\}.$  

\section{Orlik-Solomon algebras}
  
Often is it useful to realize a sequence of non-negative integers (or coefficients of a polynomial or formal power series) as the sequence of dimensions of the graded pieces of a graded algebra. The Orlik-Solomon algebra accomplishes this for the coefficients of the characteristic polynomial of a loopless matroid.  
Indeed, its Hilbert series is (after a simple algebraic transformation) the characteristic polynomial.    
The definition of the Orlik-Solomon algebra came out of the de Rham cohomology of complements of complex hyperplane arrangements.  We sketch this development to provide context and motivation.  

We shall assume a basic knowledge of the exterior algebra of a vector space. See, for example, the books 
\cite{MR1878556, MR0941522}. 

\subsection{Broken circuits}

Let $M$ be a rank-$r$ loopless matroid on the set $E$ of size $n$ with elements labeled by the integers $1,2,\ldots,n.$  
A {\em broken circuit} \index{Circuit!broken} is a subset of $E$ of the form $C - m$ where $C$ is a circuit of $M$ and $m$ is the minimum element in $C.$ A subset $I \subseteq E$ is a \emph{no-broken-circuit-} or \emph{\nbc-set} \index{nbc-set} if $I$ contains no broken circuits. 
An \nbc-set cannot contain a circuit and hence is independent.  
The collection 
of  \nbc-sets forms an abstract simplicial complex,  
the {\em broken circuit complex} of $M$ (relative to the labeling).  This complex is the subject of Chapter NN.  

In this section, we shall need two results about \nbc-sets.  The first is a handy reformulation of the definition.

\begin{lem}\label{reformulation}  A subset $i_1 i_2 \ldots i_k$ of $E,$ where $i_1 < i_2 < \cdots < i_k,$ is an  
\nbc-set if and only if, for $1 \leq t \leq k$,
$i_t$ is the minimum element in the closure $\mathrm{cl}(\{i_t, i_{t+1}, \ldots, i_k\}).$ In particular, if $M$ has rank $r,$  
then every \nbc-set of size $r$ contains $1.$     
\end{lem}

The second is a fundamental result due to Whitney \cite{MR1503085} for graphic matroids.  Recall that    
the \emph{characteristic polynomial}\index{Polynomial!characteristic}   
$\chi(M;\lambda)$ and \emph{Whitney numbers $w_k(M)$ of the second kind} are defined by 
\[
\chi (M; \lambda) = \sum_{X \in L(M)}  \mu (\emptyset,X) \lambda^{r -r(X)} = \sum_{k=0}^r (-1)^k w_k(M) \lambda^{r - k}, 
\]
where the first sum ranges over all flats $X$ in the lattice $L(M)$ of flats of $M.$   

\begin{thm}\label{face} Let $M$ be a rank-$r$ loopless matroid and $X$ be a flat of $M.$ Then 
the $(-1)^{r (X)} \mu(\emptyset,X)$ equals the number of \nbc-sets $I$ with closure equal to $X.$  In particular, $w_k(M)$ equals the number of \nbc-sets of size $k.$ 
\end{thm}

\begin{ex}\label{runningexample0}
Let $K$ be the rank-$3$ simple matroid on the set $123456$ with $3$-element circuits $123, 156, 246, 345$ shown in Figure \ref{fig:matroidK}.  The matroid $K$ is the cycle matroid of the complete graph $K_4$ on $4$ vertices and the matroid of the Coxeter hyperplane arrangement $A_3.$   Its \nbc-sets are the empty set $\emptyset,$  
all $1$-element subsets, $1, 2, 3, 4, 5, 6,$
the $2$-element subsets 
\[
12, 13, 24, 25, 34, 35, 15, 16, 14, 25, 36,
\]
and the $3$-element subsets 
\[
124, 125, 126, 134, 135, 136. 
\] 
As predicted by Theorem \ref{face},  $\chi(K ;\lambda) = \lambda^3 - 6\lambda^2 + 11 \lambda - 6. $
\end{ex}

\showhide{
\begin{figure}
  \centering
  \begin{tikzpicture}[scale=1]
  \filldraw (.3,.6) node[above] {$2$} circle  (2pt);
  \filldraw (.6,1.2) node[above] {$3$} circle  (2pt);
  \filldraw (0.9,.4) node[above] {$\,4$} circle  (2pt);
  \filldraw (0,0) node[below] {$1$} circle  (2pt);
    \filldraw (1,0) node[below] {$5$} circle  (2pt);
  \filldraw (2,0) node[below] {$6$} circle  (2pt);
  \draw[thick](0,0)--(2,0);
  \draw[thick](0,0)--(.6,1.2);   \draw[thick](1,0)--(.6,1.2);  \draw[thick](2,0)--(.3,.6);
 \end{tikzpicture}
  \caption{The matroid $K$ }
  \label{fig:matroidK}
\end{figure}
}

\subsection{Exterior and graded algebras}

Let $\k$ be a field of characteristic not equal to $2,$  $E$ be a finite set labeled by the integers $1,2,\ldots,n,$ and $\Ext (E)$ be the exterior algebra\index{Algebra!exterior} of the $|E|$-dimensional vector space $\k^E$ spanned by the standard basis 
$e_i, \, i \in E.$  
%
If $T$ is a sequence $(i_1,i_2, \ldots,i_k)$ with terms in $E,$ then the exterior product $e_T$ is defined by 
\[
e_T = e_{i_1} \wedge e_{i_2} \wedge \cdots \wedge e_{i_k}. 
\] 
Note that $e_{\emptyset} = 1.$  
The $2^n$ exterior products $e_T,$ where $T$ is an increasing sequence, form a 
a basis for $\Ext (E).$  

For $0 \leq k \leq n,$ let $\Ext^k (E)$ be the subspace (of dimension $\binom {n}{k}$) spanned by the exterior products $e_T,$ where $T$ is a length-$k$ sequence.   
If $a \in \Ext^k (E)$ for some $k,$ then $a$ is said to be {\em homogeneous.}\index{Homogeneous}  If $a \neq 0,$ it is assigned the {\em grade}\index{Grade} $k$ and we write $|a| = k.$  Under this grading, $\Ext (E)$ forms a graded algebra, in the sense that as a vector space,  
\[
\Ext (E) = \bigoplus_{k=0}^n \Ext^k (E), 
\]
and the product of a (homogeneous) element in $\Ext^j (E)$ and an element in $\Ext^k (E)$ is an element in $\Ext^{j+k} (E).$    
Multiplication in $\Ext (E)$ is {\em graded-commutative:} if $a$ and $b$ are homogeneous,  then 
\[
a \wedge b = (-1)^{|a||b|} b \wedge a. 
\] 

Let $M$ be a matroid on $E$ with lattice $L(M)$ of flats.  Then $M$ defines another grading on $\Ext (E).$  If $X$ is a flat, let $\Ext^X (E)$ be the subspace of $\Ext (E)$ spanned by the exterior products $e_T$ such that the closure $\cl (T)$ equals $X.$  Then 
\[
\Ext (E) = \bigoplus_{X \in L(M)} \Ext^X (E)
\]
and the product of an element in $\Ext^X (E)$ and an element in $\Ext^Y (E)$ is an element in $\Ext^{X \vee Y} (E),$ 
where $X \vee Y$ is the join $\cl (X \cup Y).$   If $a \in \Ext^X (E)$ for some flat $X,$ we say that $a$ is 
{\em $M$-homogeneous}\index{Homogeneous!$M$-} and has 
{\em $M$-grade}\index{Grade!$M$-}  $X.$  

%

Orlik-Solomon algebras are constructed by taking quotients of exterior algebras by an ideal\index{Ideal} defined using a boundary operator.   An ideal $I$ in $\Ext (E)$ is {\em homogeneous} (respectively, {\em $M$-homogeneous}) if $I$ is generated by homogeneous (respectively, $M$-homogeneous) elements.  The following easy proposition summarizes the underlying algebra.   

\begin{prop}\label{homogeneousgraded} Let $I$ be a homogeneous ideal in $\Ext (E).$ Then the quotient $\Ext (E) / I$ is graded, with 
\[ 
(\Ext (E)/I)^k =  \Ext (E)^k \big/  (I \cap \Ext ^k (E)). 
\]
An analogous assertion holds for $M$-homogeneous ideals. 
\end{prop}

The \emph{boundary operator} $\partial: 
\Ext (E) \to \Ext (E)$ is defined on the basis elements by $\partial 1 = 0$ and 
if $T = \{i_1,i_2, \ldots,i_k\}$ and $k > 0,$ then 
\[
\partial e_T = \sum_{j=1}^k  (-1)^{j-1}  e_{T - i_j}, 
\]   
where $T - i_j$ is the length-$(k-1)$ sequence obtained by deleting the $j$th term $i_j.$  
The definition is then extended by linearity to all of $\Ext (E).$   
We note three elementary properties.

\begin{lem}\label{elem}  
(a) $\partial\partial = 0.$  

\noindent
(b) $\partial (e_1 \wedge e_2 \wedge \cdots \wedge e_k) = (e_2 - e_1) \wedge  (e_3 - e_1) \wedge \cdots \wedge  (e_k - e_1).$  

\noindent 
(c) $\partial$ is a {\em  
graded derivation:} for homogeneous elements $a$ and $b,$   
\[
\partial(a \wedge b) = (\partial a) \wedge b + (-1)^{|a|} a \wedge (\partial b). 
\]
\end{lem}

\subsection{Orlik-Solomon algebras defined}  

We are now ready to define Orlik-Solomon algebras \cite{MR0558866}.     
  
\begin{defn}\label{OS-defn}\index{Algebra!Orlik-Solomon} 
Let $M$ be a rank-$r$ matroid on the set $E$ labeled by $1,2,\ldots,n$ and $I(M)$ be the ideal of $ \Ext(E)$ generated by the set 
\[
\{ \partial e_C \,:\, C \ \text{is a circuit of} \,  M\}.
\]
The {\em Orlik-Solomon algebra} $\OS (M)$ of $M$ is the quotient algebra defined by  
\[
\OS (M)= \Ext(E)/I(M).  
\]
We will denote by $\omega_T$ the image of $e_T$ in $\OS (M).$  
\end{defn}

\begin{ex} \label{runningexample1} 
Continuing example \ref{runningexample0}, the ideal $I(K)$ is generated by seven homogeneous elements.  Four are from the circuits $123, 156, 246, 345$ of size $3:$     
\[
e_{23} - e_{13} + e_{12},  \,\, e_{56} - e_{16} + e_{15}, \,\, 
e_{46} - e_{26} + e_{24},  \,\, e_{45} - e_{35} + e_{34}, 
\]
and three from the circuits $1245, 1346, 2356$ of size $4:$
\[
e_{245} - e_{145} + e_{125} - e_{124},\,e_{346} - e_{146} + e_{136} - e_{134},\,e_{356} - e_{256} + e_{236} - e_{235}.
\]
\end{ex}

We note some easy consequences.  If $M$ has a loop $i,$ then $\partial e_i = 1,$ $1 \in I(M),$ and $\OS (M) = 0.$  If $i$ and $j$ are parallel elements, then $\{i,j\}$ 
is a circuit, $\partial (e_i \wedge e_j) = e_i - e_j,$ and $\omega_i = \omega_j$ in $\OS (M).$  If $M$ has no loops, then $\OS (M)$ is naturally isomorphic to the Orlik-Solomon algebra of a simplification of $M.$  Thus, we may as well work with simple matroids or equivalently, as Orlik and Solomon did in \cite{MR0558866}, geometric lattices.      

\begin{lem}\label{dep}  If $D \subseteq E$ is dependent, then $e_D \in I(M)$ and $\partial e_D \in I(M).$
\end{lem}

{\it Proof.}  If $j$ is an element of a circuit $C,$ then  
\[
e_{j} \wedge \partial e_C =  e_{j} \wedge \sum_{i \in C} e_{C - i} = \pm e_C. 
\]
Hence, $e_C \in I(M).$   If $D$ is dependent, then $D$ contains a circuit $C,$ $e_D = \pm e_{D-C} \wedge e_C,$ and 
$e_D \in I(M).$   
To prove the second assertion, we use the fact that $\partial$ is a graded derivation: 
\[
\partial e_D = \partial (e_{D - C} \wedge e_C) = e_{D-C} \wedge \partial e_C + (-1)^{|D-C|} \partial e_{D-C}  \wedge e_C.      
\]
Since both terms on the right are in $I(M),$ we conclude that $\partial e_D \in I(M).$     
\hfill $\square$ \medskip \break


\begin{cor} $\OS^k(M) = 0$ if $k > r.$   
\end{cor}

It is immediate from the definition that the ideal $I(M)$ is homogeneous.  Moreover, if $C$ is a circuit and $i \in C,$ then 
$\cl (C-i) = \cl (C)$ and hence, $I(M)$ is $M$-homogeneous as well.  By Proposition \ref{homogeneousgraded},  
$\OS (M)$ is graded and $M$-graded.  
A non-zero 
homogeneous or $M$-homogeneous element has the same grade or $M$-grade in $\OS (M)$ it had in $\Ext (E).$ 
By Lemma \ref{dep}, an element $\omega_T$ is non-zero in $\OS (M)$ only if $T$ is independent.  Thus, in $\OS (M),$  
$M$-grading refines grading.   

\begin{prop}\label{Mgrade}
$\displaystyle{
\OS^k (M) = \bigoplus_{X \in L(M), \, r(X) = k}  \,\,\, \OS^X (M).  }$ 
\end{prop}

\begin{ex}\label{runningexample2}
Concluding examples \ref{runningexample0} and \ref{runningexample1}, let $X$ be the flat $123$ in $K.$  Then $\OS^X (K)$ is the $4$-dimensional vector space spanned 
by $e_{123}, e_{12},e_{13},e_{23}$ modulo the linear relations $e_{123} = 0$ and $e_{23} = e_{13} - e_{12}.$  Thus, $\OS^X (K)$ is a $2$-dimensional space with basis $e_{12}, e_{13}.$    
Now let $E =123456.$ Then, $\OS^E (K)$ is the quotient of a $20$-dimensional space (spanned by $e_{S}, S \subseteq 123456,$ $|S| = 3$) modulo 
linear relations implied algebraically by relations given in example \ref{runningexample1}.  The theory of \nbc-monomials, developed in the next subsection, will show that  
$\OS^E (K)$ has dimension $6.$    
\end{ex}

\subsection{An \nbc-basis} 

An {\em \nbc-monomial}\index{nbc-monomial} is an element in $\OS (M)$ 
of the form $\omega_J,$ where $J$ is the sequence obtained by putting an \nbc-set in increasing order.  

\begin{thm}\label{basis}  Let $M$ be a simple rank-$r$ matroid on the set $E$ labeled by $1,2,\ldots,n$ and $X$ be a flat of $M.$           Then 

\vskip 0.1in 
(a)  The \nbc-monomials $\omega_T,$ where $\cl (T) = X,$ form a basis for $\OS^X (M).$  

(b)  The \nbc-monomials $\omega_T,$ where $|T| = k,$ form a basis for $\OS^k (M).$ 

\vskip 0.1in \noindent 
In particular, the \nbc-monomials form a basis for $\OS (M).$ 
\end{thm}

{\it Proof.}  
We first show that the \nbc-monomials span by a Gr\"obner basis argument.    
 To do this, we impose the {\sf glex} or graded lexicographic order derived from the natural order on $12\ldots n$ on the set of increasing sequences $(i_1,i_2, \ldots,i_k)$ of elements from $E$ by specifying 
\[
(i_1, i_2, \ldots, i_k)  <  (j_1, j_2, \ldots, j_m)
\]
if either $k < m$ or $k = m$ and for some $t, 1 \leq t \leq k,$  $i_1 =j_1, i_2 = j_2, \ldots, i_{t-1} = j_{t-1},$ 
and $i_t < j_t$.
The {\sf glex} order can also be imposed on subsets:  $S < T$ if when put in increasing order, the sequence given by $S$ is less than the sequence given by $T.$  The extended orders are total orders with no infinite descending chains.  

By Lemma \ref{dep}, $\OS (M)$ is spanned by elements $\omega_I,$ where $I$ is an independent set.  Thus, it suffices to show that every element $\omega_I,$ where $I$ is independent, is a linear combination of \nbc-monomials. 
Let $I = \{i_1,i_2, \ldots,i_k\}$ be an independent set containing a broken circuit $C - m.$  Then $m \cup I$ is dependent, 
$ \partial e_{m \cup I} \in I(M),$ and  
\[
\omega_I = \sum_{t=1}^k (-1)^t  \omega_{ m \cup (I - i_t)}.  
\]
where by Lemma~\ref{dep}, $\omega_{m \cup (I - i_t)} \neq 0$ only if $i_t \in C$ (and $m \cup (I - i_t)$ is independent).   Since 
$m < i_t$ if $i_t \in C,$ we have written $\omega_I$ as a linear combination of elements $\omega_J,$ where $J$ is independent and $J < I.$    
Repeating this (a finite number of times), we can write $\omega_I$ as a linear combination of elements $\omega_J,$ 
where $J$ are \nbc-sets.  We have proved that the \nbc-monomials span.    

It remains to show that the \nbc-monomials are linearly independent.  By Proposition~\ref{Mgrade}, $\OS^k (M)$ decomposes into a direct sum of subspaces $\OS^X(M),$ where $r(X) = k$ and 
$\OS^X (M) =  \Ext^X (E) / (I(M) \cap \Ext^X (E)).$  
This decomposition means that a minimal linear relation among \nbc-monomials in $\OS^k (M)$ lies inside $\OS^X(M)$ for some rank-$k$ flat $X$ and is a consequence of an element in the intersection $I(M) \cap \Ext^X (E).$  
Thus, to prove assertion (b) of Theorem \ref{basis} by induction, it suffices to show assertion (a) for all rank-$k$ flats 
$X$ assuming as induction hypothesis assertion (b) for $k-1.$  
  
To begin the induction, observe that since $M$ is assumed to be simple, all circuits have at least $3$ elements.  Hence, 
$I(M) \cap \Ext^0 (E) = I(M) \cap \Ext^1 (E) = 0.$     
It follows that the $\OS^0 (M) = \Ext^0 (E)$ and $\OS^1 (M) = \Ext^1 (E)$ and the \nbc-monomials  
$\omega_{\emptyset}, \omega_1, \omega_2, \ldots, \omega_n$ are linearly independent.  

The next lemma gives a method for constructing a linear relation in $\OS^{k-1}(M)$ given one in $\OS^k (M).$  

\begin{lem} \label{boundaryonOS} 
For $1 \leq k \leq r,$ the boundary operator $\partial:\Ext^k(E) \to \Ext^{k-1} (E)$ induces a linear operator 
$\partial: \OS^k (M) \to \OS^{k-1}(M).$    
\end{lem}  

{\it Proof.} 
By Lemma \ref{elem}, $\partial \partial = 0.$  Hence, $\partial (I(M)) = 0$ in $\Ext (E)$ and $\partial$ is 
defined on the quotient $\Ext (E)/I(M).$   
\hfill $\square$ \medskip \break  
  
For the induction step, let $X$ be a rank-$k$ flat and 
\[
\sum_{S} a_S \omega_S =0
\]
be a linear relation in $\OS^X (M),$ where the subsets $S$ are \nbc-sets with closure $X.$  By Lemma 
\ref{reformulation}, 
$S = m \cup (S-m),$ where $m$ is the minimum element in $X.$ Applying $\partial,$ 
we obtain the following linear relation in $\OS^{k-1}(M):$  
\[
\sum_{S} a_S \partial \omega_S = \sum_{S} a_S \omega_{S-m} + \sum_{T} b_T \omega_T, 
\]
where on the right hand side, the first sum ranges over distinct \nbc-sets $S - m$ (not containing $m$), and the second sum 
ranges over \nbc-sets $T$ containing $m.$  
By induction, the derived linear relation is trivial and hence, $a_S = 0$ for all $S,$ that is, all linear relations among \nbc-monomials in $\OS^k (M)$ are trivial.   
\hfill $\square$ \medskip \break 

The {\em Hilbert series} $H(\R;t)$ \index{Series!Hilbert} of a graded algebra $\R = \bigoplus_{k \geq 0} \R^k$ is the formal power series 
\[
H(\R;t)=\sum_{k\geq 0} \dim (\R^k) \, t^k.
\]
Theorems \ref{face} and \ref{basis} imply the following corollary. 

\begin{cor}\label{Hilbert}   Let $M$ be a rank-$r$ simple matroid and $Y$ a flat of $M.$  Then   
\[
\dim \OS^Y (M) = (-1)^{r(Y)} \mu (\emptyset, Y) 
\]
and 
\[
H(\OS (M);t) = \sum_{X \in L(M)} \, (-1)^{r(X)} \mu (\emptyset,X) t^{r(X)} = (-t)^r \chi (M; -1/t).
\]
\end{cor}

A dual to the Orlik-Solomon algebra $\OS (M)$ is described in \cite{MR1458410}.  This algebra is constructed by taking a quotient of a vector space spanned by flags (that is, maximal chains) in    the lattice $L(M)$ of flats.  

\subsection{\OS-equivalence} 

To what extent does the algebra $\OS (M)$ determine the (simple) matroid $M$?   
The following lemma might seem to answer the question. 

\begin{lem} Let $T \subseteq E$. Then $T$ is independent in $M$ if and only if $\omega_T \neq 0$ in $\OS (M).$  
\end{lem}

{\it Proof.} If $T$ is dependent then $\omega_T = 0$ by Lemma~\ref{dep}. If $T$ is independent, then one can relabel $E$ so that $\omega_T$ is an \nbc-monomial and hence nonzero by Theorem~\ref{basis}.
\hfil $\square$ \medskip \break 
However, if $\OS (M)$ is only known as an algebra (up to isomorphism), then one may not be able to identify the elements $e_S.$ The following example, due to Rose and Terao \cite{MR1080977}, shows that this situation may occur.  

\begin{ex}\label{sixpoints}
Let $E =123456,$ $M_1$ be the simple rank-$3$ matroid on $E$ with $3$-element circuits $123$ and $456,$ and 
$M_2$ be the simple rank-$3$ matroid with $3$-element circuits $123$ and $345$ depicted in Figure \ref{fig:2guys}.  
Then $\OS (M_1)$ is isomorphic to $\OS (M_2).$

\showhide{
\begin{figure}
  \centering
  \begin{tikzpicture}[scale=1]
  \filldraw (.2,.8) node[above] {$4$} circle  (2pt);
  \filldraw (1.2,.8) node[above] {$5$} circle  (2pt);
  \filldraw (2.2,.8) node[above] {$6$} circle  (2pt);
  \filldraw (0,0) node[below] {$1$} circle  (2pt);
    \filldraw (1,0) node[below] {$2$} circle  (2pt);
  \filldraw (2,0) node[below] {$3$} circle  (2pt);
  \draw[thick](0,0)--(2,0);
  \draw[thick](.2,.8)--(2.2,.8);
  \node at (1,-0.75) {$M_1$}; 
%
  \filldraw (4.5,.6) node[above] {$4$} circle  (2pt);
  \filldraw (5,1.2) node[above] {$5$} circle  (2pt);
  \filldraw (5.7,.7) node[above] {$6$} circle  (2pt);
  \filldraw (4,0) node[below] {$1$} circle  (2pt);
    \filldraw (5,0) node[below] {$2$} circle  (2pt);
  \filldraw (6,0) node[below] {$3$} circle  (2pt);
  \draw[thick](4,0)--(6,0);
  \draw[thick](4,0)--(5.1,1.2);
  \node at (5,-0.75) {$M_2$};
 \end{tikzpicture}
  \caption{Two matroids $M_1$ and $M_2$ with isomorphic Orlik-Solomon algebras  }
  \label{fig:2guys}
\end{figure}
}
   
An explicit isomorphism $\OS (M_1) \to \OS (M_2)$ is induced by 
the isomorphism $\Phi: \Ext (E) \to \Ext (E)$ given by
\begin{eqnarray*} 
&&  e_1 \mapsto e_1, \quad   e_2 \mapsto e_2, \quad  e_3 \mapsto e_3,  \\
&&  e_4 \mapsto  e_3-e_5+e_6, \quad
e_5 \mapsto e_4-e_5+e_6, \quad   e_6 \mapsto e_6
\end{eqnarray*}
on $\Ext^1 (E).$  An easy calculation (with the help of Lemma~\ref{elem}) shows that $\Phi$ maps the ideal $I(M_1)$ into the ideal $I(M_2)$ and hence, $\Phi$ gives a homomorphism between the Orlik-Solomon algebras.  
Since the two matroids have the same characteristic polynomial
$\lambda^3 - 6\lambda^2 + 13 \lambda - 8,$  $\Phi$ is an isomorphism by Corollary~\ref{Hilbert}.
\end{ex}
 
Using the idea in example~\ref{sixpoints}, Eschenbrenner and Falk \cite{MR1719136} constructed other examples.  In particular, they show that, for any matroid $M$ and positive integer $m,$ there are $m$ non-isomorphic extensions of $M$  having isomorphic Orlik-Solomon algebras but pairwise distinct Tutte polynomials.
They also show, in the other direction, that there are matroids with the same Tutte polynomial but non-isomorphic Orlik-Solomon algebras (see Example \ref{OSPaving}).  
From these examples, one sees that the answer to the question at the start of this subsection is far from clear.

\subsection{Topology of complex hyperplane arrangements} 

The motivation behind Orlik-Solomon algebras lies in the topology of complex hyperplane arrangements. 
A hyperplane $H$ in $\mathbb{C}^r$ is the kernel of a nonzero linear form $\alpha_H,$ 
that is 
\[
H = \{(x_1,x_2, \ldots,x_r): \alpha_H (\underline{x}) = a_1x_1 + a_2x_2 +\cdots + a_rx_r = 0 \}. 
 \] 
If $\mathcal{H}$ is an {\em arrangement,}\index{Arrangement!hyperplanes} (that is, a finite set) of distinct hyperplanes, then linear dependence of the linear forms $\alpha_H, H \in \mathcal{H}$ defines a simple matroid $M$ on the set $\mathcal{H}.$  The matroid $M$ has rank $r$ if and only if the intersection $\bigcap_{H \in \mathcal{H}}  H = 0$ and we shall assume that that is the case.  
The {\em complement} $X$ of an arrangement $\mathcal{H}$ is defined by 
\[
X = \mathbb{C}^r \setminus \bigcup_{H \in \mathcal{H}} H.
\]
An important area in the study of complex hyperplane arrangements is the topology, specifically, the de Rham cohomology,  of the complement.   
Very briefly, the aim of de Rham cohomology is to calculate the {\em cohomology ring} $\mathsf{H}(X,\mathbb{C})$ constructed by taking a direct sum of cohomology groups of a cochain complex defined by differential forms.  (See, for example, \cite{MR0413152} for a detailed account.)
For complements of arrangements, $\mathsf{H}(X, \mathbb{C})$ depends only on the matroid $M$ and not on the linear forms $\alpha_H.$   
  
If $H = \ker \alpha_H,$ then the {\em logarithmic $1$-form} $\bar\omega_H$ is defined by   
\[
\bar\omega_H = d \log (\alpha_H) = \frac {d \alpha_H}{\alpha_H},
\]
where $d$ is the de Rham differential.  By the quotient rule, $d \bar\omega_H = 0$ and so $\bar\omega_H$ belongs to a cohomology class in $\mathsf{H}(X,\mathbb{C}).$  If $\mathcal{H}$ is an arrangement, then the $1$-forms $\bar\omega_H, \, H \in \mathcal{H},$ generate an algebra $\mathsf{R}(\mathcal{H})$ of differential forms with a natural map $\mathsf{R}(\mathcal{H}) \to \mathsf{H}(X,\mathbb{C}).$    
                                                                                                                      
\begin{prop} The map $\Ext (\mathcal{H}) \to \R(\mathcal{H}), \omega_H  \mapsto \bar\omega_H$ induces a surjective graded-algebra homomorphism $\OS (M) \to \R (\mathcal{H})$.
\end{prop}  

\noindent 
{\it Proof.} 
For a sequence $S= (H_1,H_2, \ldots, H_k)$ of hyperplanes in $\mathcal{H},$ let 
\[
\bar\omega_S = \bar\omega_{H_1} \wedge \bar\omega_{H_2} \wedge \cdots \wedge \bar\omega_{H_k}.  
\]
To prove the proposition, it suffices to show for every circuit $C=\{H_1, \ldots, H_k\}$ of the matroid $M$,  
$\partial \bar\omega_C $ 
is zero as a differential form on $X$. For $1 \leq t \leq k,$ we write $\alpha_{H_t} = \alpha_t$ and $\bar\omega_{H_t} = \bar\omega_t.$    
Since $\alpha_t$ can be multiplied by a nonzero complex number without changing $\bar\omega_t,$ we may assume the minimal linear relation on the circuit $C$ is 
$\alpha_{1} = \alpha_{2} + \cdots + \alpha_{k}.$  
As in Lemma~\ref{elem}, we have  
\begin{equation}
\partial \bar\omega_C = \sum_{t=1}^k (-1)^{t-1} \bar\omega_{C - H_t} = 
(\bar\omega_{2}-\bar\omega_{1}) \wedge \cdots \wedge (\bar\omega_{k} -\bar\omega_{1}).
\label{factor}
\end{equation}
Observe that 
\begin{equation}
\bar\omega_{t} - \bar\omega_{1} 
= 
d \log(\alpha_t) -d \log(\alpha_1) = d\log(f_t) = \frac{d f_t}{f_t},
\label{eqn:logs}
\end{equation}
where for $2 \leq t \leq k,$  $ f_t = \alpha_{t} / \alpha_{1}.$  Since each $f_t$ is a rational function with no poles on $X$ and  $f_2 + \cdots + f_k = 1,$  $d f_2 + \cdots + d f_k = 0.$  This dependence relation implies $df_2 \wedge \cdots \wedge df_k=0;$  with identity \eqref{eqn:logs}, this in turn implies that the product in \eqref{factor} vanishes at every point of $X.$  
\hfill $\square$  \medskip \break 

We have surjective algebra homomorphisms $\OS (M) \to \mathsf{R}(\mathcal{H}) \to \mathsf{H}(X, \mathbb{C}).$  
These homomorphisms are isomorphisms.  The isomorphism $\mathsf{R}(\mathcal{H}) \to \mathsf{H}(X, \mathbb{C})$ 
is due to Brieskorn \cite{MR0422674} and implies that $X$ is a formal space in the sense of rational homotopy theory.  The isomorphism $\OS (M) \to \mathsf{R}(\mathcal{H})$ is a deep theorem of Orlik and Solomon \cite{MR0558866} (which builds on work of Arnold \cite{MR0242196} and Breiskorn \cite{MR0422674}) which we will state without proof.    

\begin{thm} Let $\mathcal{H}$ be a complex hyperplane arrangement with matroid $M.$  Then   
$\OS (M)$ is isomorphic to $\mathsf{H}(X,\mathbb{C})$ as graded algebras. 
\end{thm}

A consequence of this theorem is that the matroid of a complex arrangement determines the cohomology ring of its complement.  
Detailed accounts of the theory of complex hyperplane arrangements can be found in \cite{Arrangementsbook, MR1217488}.  

\section{Valuative functions on matroid base polytopes} 

The Tutte polynomial defines a valuative function on polytopes defined by bases of matroids.  We give a short exposition of this area focusing on connections with the Tutte polynomial of matroids.    

\subsection{Subdivisions of base polytopes}  If $E$ is a finite set, let $\mathbb{R}^E$ be the $|E|$-dimensional real vector space with coordinates labeled by the set 
$E,$ so that the standard basis vectors $e_s,\, s \in E,$ form a basis.  
The {\em (matroid) base polytope}\index{Polytope!base} $Q(M)$ of a matroid $M$ on the set $E$ is the convex polytope in $\mathbb{R}^E$ obtained by taking the convex closure of indicator vectors of  bases of $M,$ that is,  
\[
Q(M) = \mathrm{conv}\left\{ \sum_{b \in B}  e_b:  B \,\,\text{is a basis of}\,\,M \right\}.  
\]
A {\em (base polytope) subdivision}\index{Subdivision! base polytope}  is a decomposition 
\[ Q(M) = Q(M_1) \cup Q(M_2) \cup \cdots \cup Q(M_k), 
\]  
where 
\begin{enumerate} 
\item  $Q(M_1), Q(M_2), \ldots, Q(M_k)$ are all base polytopes of matroids $M_1, M_2, \ldots, M_k$ on the set $E,$  and 
\item  for every non-empty subset $I \subseteq 12\ldots k,$ the intersection $ \bigcap_{i \in I} Q(M_i)$ is a proper face of the polytopes $Q(M_i), i \in I.$   
\end{enumerate} 

\subsection{Valuative invariants} 
A function $v$ defined on matroid base polytopes is {\em invariant}\index{Invariant}  if it depends only on the isomorphism class of the matroid.  It is 
{\em valuative}\index{Invariant!valuative}  if for every subdivision  $Q(M) = Q(M_1) \cup Q(M_2) \cup \cdots \cup Q(M_k),$  
it satisfies the inclusion-exclusion identity:  
\[
v(Q(M)) = \sum_{I \subseteq  12 \ldots k, \, I \neq \emptyset}  (-1)^{|I|-1}  v(  \bigcap_{i \in I} Q(M_i)).   
\] 
This identity is a consequence of the defining identity for valuations: for two polytopes $Q_1$ and $Q_2,$  
\[ 
v(Q_1 \cup Q_2) = v(Q_1) + v(Q_2) - v(Q_1 \cap Q_2). 
\]  
Thus, any valuation of polytopes, such as the volume or a mixed volume, gives a valuative invariant on base polytopes. 
(For a formula for the volume, see \cite{MR2610473}.)  
As a base polytope is the convex closure of indicator vectors of bases and the Tutte polynomial is a 
sum of monomials defined by internal and external activities \index{Activity} over bases, the Tutte polynomial is a valuative invariant.  
Examples of valuative functions which are not invariant can be found in \cite{MR2760656}.

\subsection{The $\mathcal{G}$-invariant}  
In \cite{MR2519849}, Derksen introduced the  $\mathcal{G}$-invariant.  It can be defined in the following way.  
Let $M$  be a rank-$r$ matroid on the set $12\ldots n.$   Let $\pi$ be a permutation on $12 \ldots n.$ The {\em rank sequence}\index{Sequence!rank}  $\underline{r}(\pi) = (r_1,r_2, \ldots,r_n)$ of $\pi$ is the sequence   
defined by $r_1 = r(\{\pi(1)\})$ and for $j \geq 2,$
\[
r_j = r(\{\pi(1), \pi(2),\ldots,\pi(j)\}) - r(\{\pi(1), \pi(2),\ldots,\pi(j-1)\}). 
\]
For matroids, $r_j = 0$ or $1,$ there are exactly $r$ $1$'s,  and the set $B(\underline{r})$ of elements $\pi(j)$ where $r_j =1$ is a basis  of $M.$   

Let $[\underline{r}]$ be a variable or formal symbol, one for each $(0,1)$-sequence $\underline{r}.$    
The {\em $\mathcal{G}$-invariant}\index{Invariant!$\mathcal{G}$-} and its coefficients $g_{\underline{r}}$ are defined by 
\[
\mathcal{G}(M) = \sum_{\pi}  [\underline{r}(\pi)] = \sum_{\underline{r}} g_{\underline{r}}(M) [\underline{r}]
\]
where the first sum ranges over all $n!$ permutations of $E.$ 
\footnote{
Derksen defines $\mathcal{G}(M)$ with $[\underline{r}]$ an element of ``a convenient basis'' of the algebra of
quasisymmetric functions.  However, the only property used is that basis elements are linearly independent and hence can be assigned values independently.  It seems simpler to define $[\underline{r}]$ as a formal symbol to which one may assign any interpretation, including a quasisymmetric function.   } 
Derksen  showed that the $\mathcal{G}$-invariant is a valuative invariant.    
Note that if $1^r0^{n-r}$ is the sequence beginning with $r$ $1$'s followed by $n-r$ $0$'s, then  
\[
g_{1^r0^{n-r}} (M) = r!(n-r)!  b(M),
\] 
where $b(M)$ is the number of bases in $M.$  

A {\em specialization} \index{Specialization!$\mathcal{G}$-} of the $\mathcal{G}$-invariant with values in an abelian group $\mathbb{A}$ is a function assigning a value in $\mathbb{A}$ to each symbol $[\underline{r}].$  
As was noted by Derksen, the formula for the rank function of the dual $M^*$ implies that the rank sequence of $\pi$ in $M^*$ can be obtained by switching $0$'s and $1$'s in the rank sequence of $\pi^{\mathrm{rev}}$ in $M,$ where $\pi^{\mathrm{rev}}$ is the permutation defined by $\pi^{\mathrm{rev}} (j) = \pi (n-j).$    Thus,  $\mathcal{G}(M^*)$ is the specialization of $\mathcal{G}(M)$ given by replacing the symbol $[\underline{r}]$ by the symbol $[\underline{r}^{*}],$ where $\underline{r}^{*}$ is the sequence obtained by reversing $\underline{r}$ and switching $0$'s and $1$'s.        

The fundamental theorem in this area, due to Derksen and Fink \cite{MR2680193}, says that the $\mathcal{G}$-invariant encompasses all valuative invariants.  
We will state this theorem without proof.  

\begin{thm}\label{DerksenFink}       
The $\mathcal{G}$-invariant is a universal valuative invariant on base polytopes, in the sense that every 
valuative invariant on base polytopes is a specialization of the $\mathcal{G}$-invariant.    In particular, the Tutte polynomial is a specialization of the $\mathcal{G}$-invariant. 
\end{thm}

Derksen \cite{MR2519849} gave an explicit specialization of the $\mathcal{G}$-invariant to the Tutte polynomial using quasisymmetric functions.  This specialization can be restated as 
\[
[r_1r_2 \ldots r_n] \mapsto \sum_{m=0}^n 
\frac {x^{r-\mathrm{wt}(r_1r_2 \ldots r_m)} y^{m-\mathrm{wt}(r_1r_2 \ldots r_m)}}{m! (n-m)!},
\] 
where the Hamming weight $\mathrm{wt}(r_1r_2 \ldots r_m)$ is the number of $1$'s in the initial segment 
$r_1r_2 \ldots r_m.$

\subsection{Paving matroids}\label{Paving} 

Compared to the Tutte polynomial, what additional information does the $\mathcal{G}$-invariant contain?
Mayhew, Newman, Welsh and Whittle \cite{MR2821559}  
conjectured that asymptotically, the proportion of  (sparse) paving matroids among all matroids tends to $1.$  Thus, if one believes the conjecture,  a ``zeroth-order'' answer can be  obtained by determining what information is contained in $\mathcal{G}$-invariants of paving matroids.  

Recall that a rank-$r$ matroid $P$ on the set $E$ is \emph{paving}\index{Matroid!paving} if all circuits have $r$ or $r+1$ elements.  A copoint $X$ (or rank-$(r-1)$ flat) is 
{\em trivial} if $|X| = r-1$ and {\em non-trivial} if $|X| \geq r.$    
A paving matroid is {\em sparse} \index{Matroid!paving!sparse} if all non-trivial copoints have size $r.$ 

Let $P$ be a paving matroid on $12 \ldots n.$   Since every subset of size $r-1$ is independent, the rank sequence of 
a permutation $\pi$ starts with 
$r-1$ $1$'s and the remaining $1$ occurs in position $i, \, i \ge r.$   If $\{\pi(1), \pi(2), \ldots,\pi(r-1)\}$ is a trivial copoint, then 
$i = r.$  If  not, then $\{\pi(1), \pi(2), \ldots, \pi(r-1)\}$ spans a non-trivial copoint $X$ and $i$ can vary from $r$ to 
$|X| + 1$ with the index equal to $i$ if and only if $\{\pi(1), \pi(2), \ldots, \pi(i-1)\} \subseteq X.$  Hence,
\begin{eqnarray*}
&&  \mathcal{G}(P) = \sum_{X \,\text{trivial}} (r-1)!(n-r+1)! [1^r0^{n-r}]
\\
&+& 
\sum_{X \,\,\text{non-trivial}} \left(
\sum_{i=r}^{|X| + 1}  \frac {|X|!}  {(|X| - i + 1)!}  (n - |X|) (n - i)! [1^{r-1 }0^{i-r}10^{n-i}]  
\right),
\end{eqnarray*}
where 
$ 1^{r-1 }0^{i-r}10^{n-i}$ is the sequence with $1$'s in the leading $r-1$ positions and the remaining $1$ in the $i$th position.   

\begin{ex}\label{examplepaving} 
Since all rank-$3$ simple matroids are paving, we have explicit formulas for their $\mathcal{G}$-invariants.  For example, let $P$ be the paving matroid on the set $123456,$ with non-trivial copoints 
$1234,456$
and trivial copoints 
$15,16,25,26,35,36.$
Geometrically, $P$ is the rank-$3$ matroid consisting of a $4$-point line $1234$ and a $3$-point line $456$ intersecting at the point $4.$  Then  
\[\mathcal{G}(P) = 48 [110010] + 132 [110100] + 540 [111000].\] 
\end{ex}

For sparse paving matroids, there is an even simpler formula derived by counting bases.   If $M$ is a rank-$r$ sparse paving matroid on $n$ elements with $\alpha$ non-trivial copoints, then $M$ has $\binom {n}{r} - \alpha$ bases.  Hence 
\[
\mathcal {G}(M) =\left(\binom {n}{r} - \alpha \right)r!(n-r)![1^r0^{n-r}] +
\alpha r!(n-r)![1^{r-1}010^{n-r-1}].
\]   
As a specific example, 
\[
\mathcal{G}(M(K_4)) = 16 \cdot 36  [111000] + 4 \cdot 36  [110100].  
\] 

The next proposition is an immediate consequence of the explicit formulas.  

\begin{prop}\label{Gpaving} The $\mathcal{G}$-invariant of a paving matroid $P$ depends only on the rank, the number of elements, and the multiset 
\[
\{ |X|: X \,\,\text{is a non-trivial copoint in}\,\,P\}. 
\]  
The $\mathcal{G}$-invariant of a sparse paving matroid depends only on the rank, the number of elements, and the number of bases.
\end{prop}

As observed in \cite{MR2768784}, the analog of Proposition \ref{Gpaving} holds for the Tutte polynomial.   
If the Mayhew-Newman-Welsh-Whittle conjecture is true, then the $\mathcal{G}$-invariant and the Tutte polynomial have the same asymptotic power to distinguish matroids.  However, the $\mathcal{G}$-invariant can distinguish pairs of matroids the Tutte polynomial cannot.  This is shown by the following example from Derksen \cite{MR2680193}.        

\begin{ex}\label{Moneandtwo}
Consider the two matroids $Q$ and $L$ in Figure \ref{fig:QL} given in \cite{MR1165543}, p.~133.  They are the smallest pair of non-paving matroids with the same Tutte polynomial.  Then 
$\mathcal{G}(Q)$ equals  
\[
48 [1010100] + 192[1011000] + 240[1100100] + 1104 [1101000] + 3456[111000],
\]
and $\mathcal{G}(L)$ equals 
\[
 24 [1010100] + 216[1011000] + 264[1100100] + 1080 [1101000] + 3456[111000].  
\]
\end{ex}

\showhide{
\begin{figure}
  \centering
  \begin{tikzpicture}[scale=1]
  \filldraw (.3,.6) node[above] {} circle  (2pt);
  \filldraw (.6,1.2) node[above] {} circle  (2pt);
  \filldraw (0.9,.4) node[above] {} circle  (2pt);
  \filldraw (0,0) node[below] {} circle  (2pt);
  \filldraw (0.7,1.25) node[below] {} circle  (2pt);
    \filldraw (1,0) node[below] {} circle  (2pt);
  \filldraw (2,0) node[below] {} circle  (2pt);
  \draw[thick](0,0)--(2,0);
  \draw[thick](0,0)--(.6,1.2);   \draw[thick](1,0)--(.6,1.2);  \draw[thick](2,0)--(.3,.6);
 \node at (1,-0.75) {$Q$}; 
  \filldraw (4.3,.6) node[above] {} circle  (2pt);
  \filldraw (4.6,1.2) node[above] {} circle  (2pt);
  \filldraw (4,0) node[below] {} circle  (2pt);
  \filldraw (4.7,1.25) node[below] {} circle  (2pt);
    \filldraw (4.8,0) node[below] {} circle  (2pt);\filldraw (5.6,0) node[below] {} circle  (2pt);
  \filldraw (6.4,0) node[below] {} circle  (2pt);
  \draw[thick](4,0)--(6.4,0);
  \draw[thick](4,0)--(4.6,1.2);   
 \node at (5.1,-0.75) {$L$};
 \end{tikzpicture}
  \caption{The matroids $Q$ and $L$ }
  \label{fig:QL}
\end{figure}
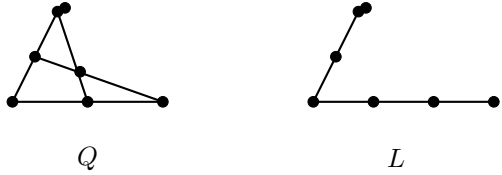
}

We end with examples from Falk \cite{MR1629681} showing that in contrast to  $\mathcal{G}$-invariants and Tutte polynomials, Orlik-Solomon algebras can distinguish some pairs of paving matroids.   

\begin{ex}\label{OSPaving}
Let $P_1$ (respectively, $P_2$) be the rank-$3$ sparse paving matroid on $123457$ with non-trivial copoints $123, 145, 356, 476$ (respectively, $123, 145, 356, 176$).  Then $\OS (P_1) \not\cong \OS (P_2).$  Falk also constructed two supersolvable rank-$3$ paving matroids with the same characteristic polynomial but non-isomorphic Orlik-Solomon algebras.  
\end{ex}
Example \ref{OSPaving} suggests the following question.  Can Orlik-Solomon algebras distinguish non-isomorphic projective planes of the same order or Dowling matroids with the same rank based on non-isomorphic groups of the same order?   

\subsection{The specialization of $\mathcal{G}$ to $T$}\label{G-Tspecialization}

By Theorem \ref{DerksenFink} or independently \cite{MR2519849}, the Tutte polynomial is a specialization of the $\mathcal{G}$-invariant.  
We discuss some consequences of this result for Tutte polynomials.     

Let $\mathcal{G}(n,r)$ be the vector space of dimension $\binom {n}{r}$ consisting of formal linear combinations of the symbols $[\underline{r}]$ where $\underline{r}$ ranges over all length-$n$ $(0,1)$-sequence with $r$ $1$'s with coefficients in $\mathbb{Q}$ (or a field of characteristic zero).  We shall construct a natural basis for this vector space.    
Let $\underline{r}$ be the length-$n$ $(0,1)$-sequence with $1$'s in positions $b_1,b_2, \ldots,b_r,$ where $b_1 < b_2 < \cdots < b_r.$ 
The {\em freedom matroid} \index{Matroid!freedom}~\footnote{Freedom matroids are first studied by Crapo \cite{MR0190045};  they have been rediscovered many times and are also known 
as nested, counting, or Schubert matroids.  Combinatorial formulas for the Tutte polynomials of freedom matroids can be found in \cite{MR2018421}. }
$F(\underline{r})$ defined by the $(0,1)$-sequence $\underline{r}$ is the rank-$r$ matroid on the set $12 \ldots n$ 
in which 
\begin{enumerate} 
\item $1,2, \ldots,b_1 - 1$ are loops (in the closure $\mathrm{cl}(\emptyset)$), 
\item for $1 \leq j \leq r-1,$   $b_j$ is added as an isthmus and $b_j, b_j + 1, \ldots,b_{j+1}-1$ are freely positioned in  
$\mathrm{cl}(\{b_1,b_2, \ldots, b_{j}\}),$ and 
\item  $b_r,b_r+1, \ldots, n$ are freely positioned in the entire matroid.  
\end{enumerate}  

We shall need a partial order on $(0,1)$-sequences.  If $\underline{r}$ and $\underline{s}$ are two length-$n$ $(0,1)$-sequences with the same number of $1$'s, then 
$\underline{s} \trianglerighteq \underline{r}$ if for every index $j,$ 
\[
s_1 + s_2 + \cdots + s_j   \geq   r_1 + r_2 + \cdots + r_j,
\] 
in other words, reading from the left, there are always at  least as many $1$'s in $\underline{s}$ as there are in $\underline{r}.$    
This order has maximum $1^r0^{n-r}$ and minimum $0^{n-r}1^r.$   
It is easy to see that the symbol $[\underline{s}]$ occurs with non-zero coefficient in $\mathcal{G}(F(\underline{r}))$ only if $\underline{s} \trianglerighteq \underline{r}.$ 
Hence, the system of equations 
\[
\sum_{\underline{s}}  g_{\underline{s}}(F(\underline{r})) [\underline{s}] = \mathcal{G}(F(\underline{r}))  
\]
is triangular with non-zero diagonal entries.  We can invert the system and write a symbol as a linear combination of $\mathcal{G}$-invariants 
of freedom matroids. This yields the following theorem.  

\begin{thm}\label{vectorspace} 
The $\mathcal{G}$-invariants $\mathcal{G}(F(\underline{r}))$ of freedom matroids form a basis for $\mathcal{G}(n,r).$  
\end{thm} 

Since the specialization of $\mathcal{G}$ to $T$ maps $\mathcal{G}(F(\underline{r}))$ to $T(F(\underline{r})),$ we have the following corollary.  

\begin{cor}\label{TPsubspace} 
The Tutte polynomial of a rank-$r$ matroid on $n$ elements is a linear combination 
of the Tutte polynomials of freedom matroids
$T(F(\underline{r})), $ 
where $\underline{r}$ ranges over all length-$n$ $(0,1)$-sequences with $r$ $1$'s.     
\end{cor}
Since almost every sufficiently large subset of a vector space spans, Corollary~\ref{TPsubspace} is almost tautological if read literally.  Its significance lies in the hope that the Tutte polynomials of freedom matroids form a natural and meaningful spanning set.  

\begin{ex}\label{TPsparsepaving}   For a rank-$r$ sparse paving matroid $M$ on $n$ elements with $\alpha$ non-trivial copoints (all of size $r$),    
\begin{eqnarray*} 
T(M) &=&  -(\alpha - 1) T(F(1^r0^{n-r})) + \alpha T(F(1^{r-1}010^{n-r-1})) 
\\
&=& -(\alpha - 1) T(U_{r,n}) + \alpha T(W^{[r]}_{r,n})
\end{eqnarray*}
where $U_{r,n}$ is a uniform matroid and $W^{[r]}_{r,n}$ is the weak-map image of $U_{r,n}$ with one non-trivial copoint with $r$ elements.  In particular,  
\[
T(M(K_4)) = -3T(U_{3,6}) + 4 T(W^{[3]}_{3,6}). 
\]
It is curious that for given $r$ and $n,$ the Tutte polynomials of the large family of sparse paving matroids lie on an affine $1$-dimensional subspace.  
\end{ex}

The assignment $\mathcal{G}(F(\underline{r}))$ to $T(F(\underline{r}))$ gives a linear transformation $\mathsf{Sp}$ from $\mathcal{G}(n,r)$ to the vector space $\mathbb{Q}[x,y]$ of polynomials in two variables $x$ and $y$ with coefficients in the rational numbers $\mathbb{Q}$ or any field of characteristic zero.    The image  is the 
subspace $\mathcal{T}(n,r)$  of $\mathbb{Q}[x,y]$ spanned by Tutte polynomials of rank-$r$ matroids on $n$ elements. Such Tutte polynomials are linear combinations of monomials $(x-1)^i(y-1)^j,$ where $0 \leq i \leq r, \, 0 \leq j \leq n-r.$  Hence an upper bound on the dimension of $\mathcal{T}(n,r)$ is $(r+1)(n-r+1).$  Thus, when $n$ is sufficiently large compared to $r,$ $\mathsf{Sp}$ is not an injection because $\dim \mathcal{G}(n,r) = \binom {n}{r}.$ 
A specific example of a {\em syzygy}\index{Syzygy} of $\mathsf{Sp},$ that is, an element in its kernel, is given by example \ref{Moneandtwo}. The linear combination  
\[
[1010100] -[1011000] - [1100100] + [1101000]  
\] 
is in the kernel because $T(M_1) = T(M_2).$  
It converts to the linear relation 
\begin{eqnarray*}
&& T(F(1010100))  - T(F(1011000)) - T(F(1100100)) 
\\
&& \qquad\qquad\qquad\qquad\qquad     + 2T(F(1101000)) - T(F(1110000)) = 0
\end{eqnarray*}
on Tutte polynomials of freedom matroids.  
A solution to the \emph{syzygy problems,} to find an explicit generating or spanning set for (1) the kernel of $\mathsf{Sp},$ and (2) the linear relations among Tutte polynomials of freedom matroids, can be found in \cite{Kungpreprint}.  This solution shows that 
$\dim \mathcal{T}(n,r) = r(n-r) + 1.$  We remark that linear relations among (arbitrary) Tutte polynomials generalize Tutte-equivalence, which are linear relations between two Tutte polynomials.

\subsection{The $F$-invariant}

A precursor of the $\mathcal{G}$-invariant is the $F$-invariant defined by Billera, Jia, and Reiner \cite{MR2552658}.  Our description of the $F$-invariant assumes an acquaintance with the greedy algorithm axiomatization (Edmonds \cite{MR0297357};  see, for example,   \cite[p.~55]{MR2849819}).  
 
Let $M$ be a matroid on $E.$  A function $f: E \to \{1,2,3,\ldots\}$ is \emph{$M$-generic} \index{Function!$M$-generic} if the minimum 
$$
\min \left\{\sum_{b \in B} f(b) : B \,\,\text{is a basis of}\,\,M  \right\}
$$
is achieved by exactly one basis.  Let $x_1,x_2, x_3,  \ldots  \,$ be variables and define the 
{\em $F$-invariant}\index{Invariant!$F$-} $F(M)$ to be formal 
power series given by  
\[
F(M) = \sum_{f} \prod_{i \in E}  x_{f(i)}, 
\]
where the sum ranges over all $M$-generic functions. 
The power series $F(M)$ is a quasisymmetric function.  As a sum over bases, it is a valuative invariant.     
Derksen has described the specialization of $\mathcal{G}$ to $F$ \cite{MR2519849}.   

\section{Coalgebras associated with matroids} 

The Tutte polynomial and other invariants satisfy identities which are manifestations of an underlying coalgebra structure.          
We give a brief and informal account of coalgebras constructed from matroids.      

Coalgebras \index{Coalgebras} are opposites of algebras.  Instead of a (bilinear) multiplication $\mathsf{R} \otimes \mathsf{R} \to \mathsf{R}, \, (x,y) \mapsto x \circ y$ in an algebra $\mathsf{R}$ ``merging'' two elements into one, a {\em combinatorial coalgebra} $\mathsf{C}$ is defined by a \emph{comultiplication} 
\[
\Delta: \mathsf{C} \to \mathsf{C} \otimes \mathsf{C} \quad    x \mapsto \sum   
\,\,  x^{\prime} \otimes x^{\prime\prime},  
\]
where the sum  ranges over ``all'' decompositions of $x$ into two parts $x^{\prime}$ and $x^{\prime\prime}.$  
%
Combinatorial coalgebras may have additional structure:  A {\em bialgebra} is a coalgebra with a multiplication $\circ$ compatible with the comultiplication and a {\em Hopf algebra} is a bialgebra with an ``antipode'', a comultiplication analog of the inverse.        
See \cite{MR3097867} for a quick introduction with detailed definitions.  

The \emph{ restriction-contraction coalgebra} \index{Coalgebras!restriction-contraction} of a minor-closed 
class $\mathcal{C}$ of matroids is defined by the comultiplication  
\[
\Delta \, M = \sum_{A \subseteq E}   \,\,  M|A  \otimes  M/A.   
\]
on the vector space  of formal linear combinations of isomorphism-classes of matroids in $\mathcal{C}.$   Usually, one imposes a multiplication $\circ,$ defined by  
\[
M \circ N = M \oplus N,
\]
to form a bialgebra.  By a method of Schmitt \cite{MR1303288}, restriction-contraction bialgebras have antipodes and are Hopf algebras.  Crapo and Schmitt \cite{MR2155892} have shown that the restriction-contraction bialgebra  of freedom matroids \index{Matroid!freedom} is generated freely as a bialgebra by a loop and a coloop.    

Many natural invariants of matroids are \emph{compatible} with the restriction-contraction coalgebra.  For example, the corank-nullity polynomial of a rank-$r$ matroid $M$ on the set $E$ satisfies the convolution-multiplication identity 
\cite{MR2718681}
\[
R(M;xy, \lambda \xi) = \sum_{A \subseteq E}  
\lambda^{r - r(A)} (-y)^{|A| - r(A)}  R(M|A;-x,-\lambda)  R(M/A;y,\xi).  
\] 
This identity converts multiplication of variables into a convolution and is a typical example of compatibility.  The corank-subset or multivariate Tutte polynomial, the $\mathcal{G}$-invariant, and the $F$-invariant satisfy similar identities.  Intuitively, the reason is that these invariants can be expressed as sums over subsets or chains of subsets and such sums can be decomposed according to the restriction-contraction comultiplication.    

The Hopf algebra $\mathsf{QSym}$ of quasisymmetric functions (see, for example, \cite{MR3097867}) is used in a fundamental way  the study of the $F$- and $\mathcal{G}$-invariants.  Indeed, $\mathcal{G}$ defines a Hopf algebra homomorphism 
from the restriction-contraction coalgebra on all matroids into $\mathsf{QSym}$ \cite{MR2519849}.  In addition, 
from the $F$-invariant, Luoto has constructed a ``matroid-friendly'' basis for $\mathsf{QSym}$ \cite{MR2417021}.   

Coalgebra-compatibility may involve other binary operations and coalgebras.  Let $P(G;x)$ be the chromatic polynomial 
\index{Polynomial!chromatic} of the graph $G$ with vertex set $V.$  If $U \subseteq V,$ let $G|U$ be the \emph{induced} subgraph on the vertex set $U$ (with all the edges in $G$ having both endpoints in $U$).   Tutte \cite{MR0223272} observed that 
\[
P(G;x+y) = \sum_{U \subseteq V}  P(G|U;x)P(G|(V \backslash U); y).
\]
Thus, the chromatic polynomial (under addition of variables) is compatible with the ``Boolean coalgebra'' defined by the comultiplication  
\[
\Delta V = \sum_{U \subseteq V}  U \otimes (V \backslash U).  
\]

\vskip 1in \noindent 
{\bf Acknowledgment.}  The authors thank Joseph Bonin, Graham Denham, Alex Fink, and Vic Reiner for many comments and suggestions.


\printindex
\end{document}